\documentstyle[12pt]{article}                  % Use this for Latex (2.09)
\input{amssymb.sty}
\makeatletter
\def\nothing#1{}
\newdimen\earraycolsep
\renewcommand{\thetable}{\arabic{table}}
\renewcommand{\thefigure}{\arabic{figure}}
\renewcommand{\title}[1]{%
  \vspace*{120\p@}%
  {\parindent \z@ \raggedright \reset@font
    \bfseries #1\par
    \nobreak
    \vskip 36\p@
  }}
\def\author#1{{\pretolerance=10000 \raggedright \advance \leftskip by 1in
\noindent #1 \vskip 1pc}}
\def\affiliation#1{{\advance\leftskip by 1in \noindent #1 \vskip -1pc}}
\def\refnote#1{{$^{\hbox{\scriptsize #1}}$}}
\def\affnote#1{\llap{$^{\hbox{\scriptsize #1}}$}}

\renewcommand\section{\@startsection{section}{1}{\z@}{2pc \@plus
      1ex minus .2ex}{1pc \@plus .2ex}{\reset@font
      \normalsize\bfseries\noindent
      {\addtocounter{section}{1}}\arabic{section}\
      {\setcounter{subsection}{0}
      \setcounter{subsubsection}{0}\setcounter{equation}{0}} }}
\renewcommand\subsection{\@startsection{subsection}{2}{\z@}{1pc \@plus 1ex
    minus.2ex}{1pc \@plus .2ex}
    {\reset@font\normalsize\bfseries
    \noindent{\addtocounter{subsection}{1}}%
    {\setcounter{subsubsection}{0}}\arabic{section}.\arabic{subsection}\ }}
\renewcommand\subsubsection{\@startsection{subsubsection}{3}{\parindent}
        {1pc \@plus 1ex minus.2ex}{-0.5em}{\reset@font\normalsize\bfseries%
        {\addtocounter{subsubsection}{1}} \hspace*{.6cm}
        \arabic{section}.\arabic{subsection}.\arabic{subsubsection}
        \hspace*{-7mm}}}

\def\AmS{{\protect\the\textfont2%
        A\kern-.1667em\lower.5ex\hbox{M}\kern-.125emS}}

\def\p@LaTeX{{\family{times}\series{m}\shape{n}\selectfont
L\kern-.36em\raise.3ex\hbox{\scriptsize A}\kern-.15em
T\kern-.1667em\lower.7ex\hbox{E}\kern-.125emX}}

\newlength{\colwidth}

\def\@oddhead{\hfil}
\def\@evenhead{\hfil}
\def\@oddfoot{{\bfseries\hfil\thepage}}
\def\@evenfoot{{\bfseries\thepage\hfil}}
\def\fnum@figure{\footnotesize\raggedright{\bfseries \figurename~\thefigure.}}
\def\fnum@table{\normalsize\raggedright{\bfseries \tablename~\thetable.}}
\long\def\@makecaption#1#2{\vskip 10\p@ {#1 #2\par}}
\long\def\@makefntext#1{\setbox0=\hbox{$\m@th^{\@thefnmark}$}\noindent
\hangindent=\wd0 \box0 #1}
\newbox\@atbox
\long\def\atable#1#2#3{\begin{table}[tbp]\centering\footnotesize
\setbox\@atbox\hbox{#2}
\parbox{\wd\@atbox}{\caption{#1}}\par\smallskip
#2
\par\smallskip\parbox{\wd\@atbox}{\raggedright #3}
\end{table}}

\newtheorem{theorem}{Theorem}

\newtheorem{proposition}[theorem]{Proposition}

\newtheorem{definition}[theorem]{Definition}

\newcommand{\ra}{\rightarrow}
\newcommand{\fl}{\forall}
\newcommand{\wt}{\widetilde}

\newcommand{\s}{\sigma}
\newcommand{\D}{\Delta}

\newcommand{\ot}{\otimes}

\newcommand{\Hc}{\mathcal{H}}
\newcommand{\g}{\gamma}

\newcommand{\ve}{\varepsilon}
\newcommand{\Cb}{\mathbb{C}}

\def\Cb{{\mathbb C}}

\def\Hc{{\cal H}}

\def\g{\gamma}

\def\s{\sigma}

\def\ve{\varepsilon}

\def\D{\Delta}

\def\fl{\forall}

\def\ot{\otimes}
\def\part{\partial}

\def\ra{\rightarrow}

\def\text{\hbox}

\def\boxit#1#2{\setbox1=\hbox{\kern#1{#2}\kern#1}%
\dimen1=\ht1 \advance\dimen1 by #1 \dimen2=\dp1 \advance\dimen2 by #1
\setbox1=\hbox{\vrule height\dimen1 depth\dimen2\box1\vrule}%
\setbox1=\vbox{\hrule\box1\hrule}%
\advance\dimen1 by .4pt \ht1=\dimen1
\advance\dimen2 by .4pt \dp1=\dimen2 \box1\relax}

\def\@nbibitem#1{\noindent \hangindent=2pc \hangafter=1
\refstepcounter{enumi}\hbox to 2pc{\arabic{enumi}.\hfil}%
\immediate\write\@auxout{\string\bibcite{#1}{\arabic{enumi}}}}
\def\numbibliography{%
\section*{REFERENCES}%
\bgroup\footnotesize
\setcounter{enumi}{0}%
\def\newblock{\hskip .11em plus.33em minus.07em}%
\let\bibitem\@nbibitem}
\def\endnumbibliography{\par\egroup}
\makeatother
\begin{document}

\begin{center}
{ \bf CYCLIC COHOMOLOGY, HOPF ALGEBRAS  \\ AND THE MODULAR THEORY}
\end{center}

\vspace{1cm}

\author{\bf Alain CONNES\refnote{1} and Henri MOSCOVICI\refnote{2}}

\affiliation{\affnote{1}  Coll\`ege de France,
3, rue Ulm,
75005 PARIS\\
and\\
$^2$Department of Mathematics, The Ohio State University\\
231 W 18th Avenue, Columbus, OH 43210 USA                     
}

\vspace{1cm}
\begin{abstract}

We associate canonically a cyclic module to 
any Hopf algebra 
endowed with a modular pair,
consisting of a group-like element and
a character, in involution. This provides the key construct allowing to
extend cyclic cohomology to Hopf algebras in the non-unimodular
case and further to develop a theory of
characteristic classes for actions of Hopf algebras compatible
not only with traces but also with the modular theory of
weights.
It applies to ribbon and to coribbon algebras, as well as 
to quantum groups and their duals.
 
\end{abstract}
\vspace{1cm}

\section{Introduction}

We solve in this paper the question left open in 
\cite{CM} (and \cite{CM2})
of extending cyclic cohomology to Hopf algebras
 in the general non-unimodular
case. The setup employed in \cite{CM} was still
 based on a partial unimodularity
condition, which was breaking the natural 
symmetry between a Hopf algebra 
and its dual.

Our solution relies on the modular theory of 
weights instead of traces.
We use this theory in an algebraic way, 
by introducing the notion of a 
$\sigma$-trace on an algebra on which our 
Hopf algebra acts, where $\sigma$
is a group-like element.
This leads to a natural construction of a cyclic module
associated to any
Hopf algebra endowed with a {\it modular pair in involution},
i.e. with a group-like element and a character such that the
corresponding doubly twisted antipode has square the identity.
The simplicial structure subjacent to this cyclic module involves the
coproduct of the Hopf algebra and the group-like element, while the
cyclic structure makes use of
the product and of the twisted antipode.

The tracial case of this construction, 
which corresponds to the group-like element being trivial,
was initially introduced in \cite{CM} under a restrictive assumption,
then  amended in \cite{CM2}; cf. also
\cite{Cr}, where it was recast in the
Cuntz-Quillen formalism -- the non-unimodular case can likewise be 
reformulated.
By extending the construction of the cyclic module to the general
non-unimodular case, we settle a problem left open in the above mentioned 
papers,
thus laying the groundwork for a theory of characteristic classes
for actions of Hopf algebras compatible with the modular
theory of weights.

It is important to mention that the non-unimodular case does arise
even in the simplest examples and that neglecting its (often hidden) presence
would give rise to misleading answers.   
The property of the existence of a modular pair in involution
is intrinsically satisfied by both the ribbon and the coribbon Hopf algebras,
as well as by the quantum groups and their duals.

\section{Characteristic classes for actions of Hopf algebras}

In what follows 
we fix a {\it modular pair}, consisting of a group-like element $\sigma$ and a
character $\delta$ of $\Hc$ such that 
$$\delta(\sigma)=1.$$ 
They 
will play the role of the module of locally compact groups.

We then introduce the twisted antipode,
\begin{equation}\label{bc}
\wt S (y) = \sum  \delta (y_{(1)})  S (y_{(2)}) \ , \ y \in {\Hc}  , \
\D  
y = \sum  y_{(1)} \ot y_{(2)}.
\end{equation}

 Given an algebra $A$, an 
action of the Hopf algebra $\Hc$ on $A$ is given by a linear map,
$$
\Hc \ot A \ra A, \quad h \ot a \ra h(a) 
$$
satisfying $h_1 (h_2  a) = (h_1  h_2) (a)$, $\fl  h_i \in {\Hc}$,
$a \in A$ and
\begin{equation}\label{ba}
h(ab) = \sum  h_{(1)}  (a)  h_{(2)}  (b)  \qquad \fl  a,b  \in A,
h \in {\Hc}.
\end{equation}
where the coproduct of $h$ is,
\begin{equation}\label{bb}
\D(h)=  \sum  h_{(1)}  \ot  h_{(2)} 
\end{equation}

\begin{definition}
We shall say that a linear form $\tau$ on $A$ is a $\sigma$-trace
under the action of $\Hc$ iff one has,
$$
\tau (ab) = \tau (b \sigma (a)) \qquad \fl  a,b  \in A  .
$$

We shall say that a $\sigma$-trace $\tau$ on $A$ is $\delta$-invariant
under the action of $\Hc$ iff 
$$
\tau (h(a)b) = \tau (a  \wt S (h)(b)) \qquad \fl  a,b  \in A  , \ h
\in 
{\Hc}. 
$$
\end{definition}
As in \cite{CM} the definition of the cyclic 
complex for $HC^*_{(\delta,\sigma)} ({\Hc})$ is uniquely 
dictated in such a way that
the following proposition holds,

\begin{proposition}
Let $\tau$ be a
$\delta$-invariant $\sigma$-trace 
on $A$, then the following defines a canonical map from $HC^*_{(\delta,\sigma)}
({\Hc})$ to 
$HC^* (A)$,
$$
\matrix{
\g (h^1 \ot \ldots \ot h^n) \in C^n (A)  , \ \g (h^1 \ot \ldots \ot
h^n) 
(x^0 , \ldots , x^n) = \cr
\cr
\tau (x^0  h^1 (x^1) \ldots h^n (x^n)). \cr
}
$$
\end{proposition}

We shall show below that the required cyclic complex can be implemented
whenever the modular pair 
${(\delta,\sigma)}$ satisfy a natural involutive condition.

\section{The cyclic module of a Hopf algebra}

In this section we shall associate a cyclic complex (in fact a 
$\Lambda $-module, where $ \Lambda$ is the cyclic category), 
to any Hopf algebra $\Hc$ (over $\Cb$) endowed with a {\it modular pair}
${(\delta,\sigma)}$ {\it in involution}, i.e. satisfying 
\begin{equation}\label{ch}
(\sigma^{-1} \wt{S})^2 = I.
\end{equation}
With the standard notation for unit $\eta : \Cb 
\ra \Hc$, counit $\ve : \Hc \ra \Cb$ and antipode $S : \Hc \ra \Hc$, 
we recall that the $\delta$-twisted antipode was defined as
\begin{equation}\label{ca}
\wt S (h) = \sum_{(h)} \delta (h_{(1)}) \ S (h_{(2)}) \quad , 
\quad h \in \Hc.
\end{equation}

\noindent The elementary properties of $S$ imply immediately that $\wt S$ is an 
algebra antihomomorphism
\begin{equation}\label{cb}
\matrix{
&\wt S (h^1  h^2) = \wt S (h^2)  \wt S (h^1) \quad , \quad \fl  
h^1 , h^2 \in \Hc \cr \cr
&\wt S (1) = 1, \hfill \cr
}
\end{equation}
a coalgebra twisted antimorphism
\begin{equation}\label{cc}
\D  \wt S (h) = \sum_{(h)} S (h_{(2)}) \ot \wt S (h_{(1)}) \quad , 
\quad \fl  h \in \Hc;
\end{equation}
and also that it satisfies
\begin{equation}\label{cd}
\ve \circ \wt S = \delta.
\end{equation}

By transposing and twisting by $\sigma$ the standard 
simplicial operators underlying the 
 Hochschild homology complex of an algebra,
one associates to $\Hc$, viewed only as a coalgebra, the 
following cosimplicial module
$\{ \Hc^{\ot n} \}_{n \geq 1}$, with face operators $\delta_i: 
\Hc^{\ot n-1} \ra \Hc^{\ot n}$,
\begin{eqnarray}\label{ce}
&&\delta_0 (h^1 \ot \ldots \ot h^{n-1}) = 1 \ot h^1 
\ot \ldots \ot h^{n-1} \nonumber \\
&& \nonumber \\
&&\delta_j (h^1 \ot \ldots \ot h^{n-1}) = h^1 \ot \ldots \ot \D h^j \ot 
\ldots \ot h^n,\  \fl  1 \leq j \leq n-1, \nonumber\\
&&  \\
&&\delta_n (h^1 \ot \ldots \ot h^{n-1}) = h^1 \ot \ldots \ot h^{n-1}
\ot \sigma \nonumber
\end{eqnarray}
and degeneracy operators $\s_i : \Hc^{\ot n+1} \ra \Hc^{\ot n}$,
\begin{equation}\label{cf}
\s_i (h^1 \ot \ldots \ot h^{n+1}) = h^1 \ot \ldots \ot \ve (h^{i+1}) 
\ot \ldots \ot h^{n+1} \ , \ 0 \leq i \leq n.
\end{equation}
The remaining two essential features of a Hopf algebra 
-- \textit{product} and \textit{antipode} -- are now brought into play, to
define the 
\textit{cyclic operators} $\tau_n : \Hc^{\ot n} \ra \Hc^{\ot n}$,
\begin{equation}\label{cg}
\tau_n (h^1 \ot \ldots \ot h^n) = (\D^{n-1}  \wt S (h^1)) \cdot h^2 
\ot \ldots \ot h^n \ot \sigma.
\end{equation}

\begin{theorem}
Let $\Hc$ be a Hopf algebra endowed
with a modular pair
${(\delta,\sigma)}$ in involution (\ref{ch}).
Then $\Hc_{(\delta,\sigma)}^{\natural} = \{ \Hc^{\ot n} \}_{n \geq 1}$ equipped 
with the operators given by (\ref{ce})--(\ref{cg}) defines a module over the
cyclic category $\Lambda$.
\end{theorem}
\smallskip
\noindent{\textit{Proof.}} The simplicial relations are easy to check
and follow from the group-like property of $\sigma$.
We shall verify, following closely the corresponding computations in \cite{CM2},
the remaining relations of the cyclic category:
\begin{eqnarray}
&&\tau_n  \delta_i = \delta_{i-1}  \tau_{n-1} \ , \ 1 \leq i \leq n  , 
\nonumber\\
&& \label{ci}\\
&&\tau_n  \delta_0 = \delta_n  ,\nonumber
\end{eqnarray}
\begin{eqnarray}
&&\tau_n  \s_i = \s_{i-1}  \tau_{n+1} \ , \ 1 \leq i \leq n, \nonumber\\
&& \label{cj} \\
&&\tau_n  \s_0 = \s_n  \tau_{n+1}^2  , \nonumber\\
&& \nonumber\\
&&\tau_n^{n+1} = I_n.\label{ck}
\end{eqnarray}

\noindent As in \cite{CM2}, we shall only use the basic properties of the 
product, the coproduct, the antipode and of the twisted 
antipode (cf. (\ref{ca})--(\ref{cd})), and adhere to 
the standard notational 
conventions for the Hopf algebra calculus (cf. \cite{S}).

We first look at the case $n=2$. Thus,
\begin{eqnarray*}
\tau_2 (h^1 \ot h^2)& = & \ \D  \wt S (h^1) \cdot h^2 \ot \sigma = \cr
&= & \ \sum \wt S (h^1)_{(1)}  h^2 \ot \wt S (h^1)_{(2)} \sigma \cr
&= & \ \sum S (h_{(2)}^1)  h^2 \ot \wt S (h_{(1)}^1) \sigma  .
\end{eqnarray*}
Its square is therefore:
\begin{eqnarray*}
\tau_2^2 (h^1 &\ot& h^2) =  \ \sum S (S (h_{(2)}^1)_{(2)}  
h_{(2)}^2)  \wt S (h_{(1)}^1) \sigma \ot \wt S (S (h_{(2)}^1)_{(1)}  
h_{(1)}^2) \sigma \cr
&= & \ \sum S (S (h_{(2)(1)}^1)  h_{(2)}^2)  \wt S (h_{(1)}^1) \sigma \ot 
\wt S (S (h_{(2)(2)}^1)  h_{(1)}^2) \sigma \cr
&= & \ \sum S (h_{(2)}^2)  (S \circ S)  (h_{(2)(1)}^1)  \wt S 
(h_{(1)}^1) \sigma  \ot \wt S (h_{(1)}^2)  (\wt S \circ S)  
(h_{(2)(2)}^1) \sigma  \cr
&= & \ \sum S (h_{(2)}^2) \ \hbox{\boxit{6pt}{$S(S(h_{(1)(2)}^1))  
\wt S (h_{(1)(1)}^1)$}} \ \sigma \ot \wt S (h_{(1)}^2)  \wt S (S 
(h_{(2)}^1)) \sigma  .
\end{eqnarray*}
The term in the box is computed as follows. With $k = h_{(1)}^1$, one 
has
\begin{eqnarray*}
\sum  S (S(k_{(2)}))  \wt S (k_{(1)})& = & \ \sum S (S(k_{(2)})) 
 \delta (k_{(1)(1)})  S (k_{(1)(2)}) \cr
&= & \ \sum S(S(k_{(2)(2)})  \delta (k_{(1)})  S (k_{(2)(1)}) = \cr
&= & \ \sum \delta (k_{(1)})  S \left( \sum k_{(2)(1)}  S (k_{(2)(2)}) 
\right) \cr
&= & \ \sum \delta (k_{(1)})  S (\ve (k_{(2)})  1) = \cr
&= & \ \sum \delta (k_{(1)})  \ve (k_{(2)}) = \delta \left( \sum
k_{(1)}  
\ve (k_{(2)}) \right) \cr
&= & \ \delta (k)  .
\end{eqnarray*}
It follows that
\begin{eqnarray*}
\tau_2^2 (h^1 \ot h^2)& = & \ \sum S (h_{(2)}^2) \sigma   \underbrace{\delta 
(h_{(1)}^1) \ot \wt S (h_{(1)}^2)  \wt S (S (h_{(2)}^1))} \sigma  \cr
&= & \ \sum S (h_{(2)}^2) \sigma \ot \wt S (h_{(1)}^2)  \wt S (\wt S (h^1)) \sigma =
\cr
&= & \ \sum S (h_{(2)}^2) \sigma  \ot \wt S (h_{(1)}^2) \wt S^2 (h^1) \sigma  ,
\end{eqnarray*}
 Thus
\begin{eqnarray*}
\tau_2^2 (h^1 \ot h^2) &= & \ \sum S (h_{(2)}^2) \ot \wt S (h_{(1)}^2) 
\cdot \sigma \ot \wt S^2 (h^1) \sigma \cr
&= & \  \D  \wt S (h^2) \cdot \sigma \ot \wt S^2 (h^1) \sigma   . 
\end{eqnarray*}
In a similar fashion,
\begin{eqnarray*}
\tau_2^3 (h^1 \ot h^2) &= & \ \sum S (S(h_{(2)}^2 )_{(2)} \sigma)  \wt S 
(h_{(1)}^2)  \wt S^2 (h^1) \sigma \ot \wt S (S (h_{(2)}^2)_{(1)} \sigma ) \sigma \cr
&= & \ \sum S (S(h_{(2)(1)}^2) \sigma )  \wt S (h_{(1)}^2)  \wt S^2 (h^1)\sigma  \ot \wt S 
(S (h_{(2)(2)}^2) \sigma ) \sigma \cr
&= & \ \sum S (S(h_{(2)(1)}^2 ) \sigma )  \wt S (h_{(1)(1)}^2)  \wt S^2 (h^1) \sigma \ot \sigma^{-1} \wt S (S (h_{(2)}^2)) \sigma \cr
&= & \ \sum \delta (h_{(1)}^2) \sigma^{-1} \wt S^2 (h^1) \sigma \ot  
\sigma^{-1} \wt S (S (h_{(2)}^2)) \sigma = \cr
&= & \  \sigma^{-1} \wt S^2 (h^1) \sigma \ot \sigma^{-1} \wt{S}^2 (h^2) \sigma 
 = h^1 \ot h^2  .
\end{eqnarray*}

We now pass to the general case. With the standard conventions of 
notation,

\begin{eqnarray*}
&&\tau_n (h^1 \ot h^2 \ot \ldots \ot h^n) = \D^{(n-1)}  \wt S (h^1) 
\cdot h^2 \ot \ldots \ot h^n \ot \sigma \cr
&&\quad = \sum S (h_{(n)}^1)  h^2 \ot S (h_{(n-1)}^1)  h^3 \ot \ldots \ot 
S (h_{(2)}^1)  h^n \ot \wt S (h_{(1)}^1) \sigma  .
\end{eqnarray*}
Upon iterating once
\begin{eqnarray*}
\tau_n^2 (h^1 \ot \ldots \ot h^n)& = & \ \sum S (S(h_{(n)}^1)_{(n)}  
h_{(n)}^2 ))  S (h_{(n-1)}^1)  h^3 \ot \cr
&\ot & \ S (S (h_{(n)}^1)_{(n-1)}  h_{(n-1)}^2))  S (h_{(n-2)}^1) 
 h^4 \ot \ldots \cr
\ldots& \ot & \ S (S (h_{(n)}^1)_{(2)}  h_{(2)}^2))  \wt S 
(h_{(1)}^1) \sigma \ot \wt S (S (h_{(n)}^1)_{(1)}  h_{(1)}^2) \sigma \cr
&= & \ \sum S (h_{(n)}^2)  S (S (h_{(n)(1)}^1))  S (h_{(n-1)}^1) 
 h^3 \ot \cr
&\ot & \ S (h_{(n-1)}^2)  S (S (h_{(n)(2)}^1))  S (h_{(n-2)}^1)  
h^4 \ot \ldots \cr
\ldots &\ot & \ S (h_{(2)}^2)  S (S (h_{(n)(n-1)}^1))  \wt S 
(h_{(1)}^1) \sigma \ot \cr
&\ot & \ \wt S (h_{(1)}^2)  \wt S (S(h_{(n)(n)}^1)) \sigma = \cr
&= & \ \sum S (h_{(n)}^2)  S (h_{(n-1)}^1  S (h_{(n)}^1))  h^3 
\ot \cr
&\ot & \ S (h_{(n-1)}^2)  S (h_{(n-2)}^1  S (h_{(n+1)}^1))  h^4 
\ot \ldots \cr
\ldots &\ot & \ S (h_{(2)}^2)  S (S (h_{(2n-2)}^1 )) \cdot \wt S 
(h_{(1)}^1) \sigma \ot \cr
&\ot & \ \wt S (h_{(1)}^2)  \wt S ( S (h_{(2n-1)}^1 )) \sigma  . 
\end{eqnarray*}
We pause to note that
$$
\sum h_{(n-1)}^1  S (h_{(n)}^1) = \sum h_{(n-1)(1)}^1  S 
(h_{(n-1)(2)}^1)
$$
equals
$$
\ve  (h_{(n-1)}^1)  1  ,
$$
after resetting the indexation. Next
$$
\sum \ve (h_{(n-1)}^1)  h_{(n-2)}^1
$$
gives $h_{(n-2)}^1$ after another resetting. In turn
$$
\sum h_{(n-2)}^1  S (h_{(n-1)}^1) 
$$
equals
$$
\ve (h_{(n-2)}^1)  1  ,
$$
and the process continues.

In the last step,

\begin{eqnarray*}
& \ &\sum S (h_{(n)}^2)  h^3 \ot S (h_{(n-1)}^2)  h^4 \ot \ldots 
\cr
\ldots &\ot & \ S (h_{(2)}^2) \ \hbox{\boxit{6pt}{$S(S(h_{(2)}^1))
\delta 
(h_{(1)(1)}^1)  S (h_{(1)(2)}^1)$}} \ \sigma \ot \cr
&\ot & \ \wt S (h_{(1)}^2)  \wt S ( S (h_{(3)}^1)) \sigma \cr
&= & \ \sum S (h_{(n)}^2)  h^3 \ot S (h_{(n-1)}^2)  h^4 \ot \ldots
\ot S 
(h_{(2)}^2) \sigma \ot \wt S (h_{(1)}^2) \wt S^2 (h^1) \sigma \cr
&= & \ \sum S (h_{(n)}^2) \ot  S (h_{(n-1)}^2) \ot \ldots \ot \wt S
(h_{(1)}^2)
\cdot 
h^3 \ot h^4 \ot \ldots \ot \sigma \ot \wt S^2 (h^1) \sigma \cr
&= & \ \D^{(n-1)} \wt S (h^2) \cdot h^3 \ot h^4 \ot \ldots \ot \sigma \ot 
\wt S^2 (h^1) \sigma  ,
\end{eqnarray*}
with the boxed term simplified as before. 

By induction, one obtains 
for any $j = 1, \ldots , n+1$,
$$
\tau_n^j (h^1 \ot \ldots \ot h^n) = \D^{n-1} \wt S (h^j) \cdot h^{j+1} 
\ot \ldots \ot h^n \ot \sigma \ot \ldots \ot \wt S^2 (h^{j-1}) \sigma,
$$
in particular
$$
\tau_n^{n+1} (h^1 \ot \ldots \ot h^n) = \D^{n-1} \wt S (\sigma) \cdot \wt S^2 
(h^1) \sigma \ot 
\ldots \ot \wt S^2 (h^n) \sigma = h^1 \ot \ldots \ot h^n  .
$$

The verification of the compatibility relations (\ref{ci}), (\ref{cj}) is 
straightforward. Indeed, starting with the compatibility
with the face operators, one has:
\begin{eqnarray*}
\tau_n  \delta_0 (1 \ot h^1 \ot \ldots \ot h^{n-1})& = & \ \tau_n (1
\ot 
h^1 \ot \ldots \ot h^{n-1}) = \cr
&= & \ \D^{n-1}  \wt S (1) \cdot h^1 \ot \ldots \ot h^{n-1} \ot \sigma 
\cr
&= & \ h^1 \ot \ldots \ot h^{n-1} \ot \sigma \cr
&= & \ \delta_n (h^1 \ot \ldots \ot h^{n-1})  , 
\end{eqnarray*}
then
\begin{eqnarray*}
\tau_n  \delta_1 (h^1 \ot \ldots &\ot& h^{n-1}) =  \ \tau_n \,( \D  h^1 
\ot h^2 \ot \ldots \ot h^{n-1}) \cr
&= & \ \sum \tau_n (h_{(1)}^1 \ot h_{(2)}^1 \ot h^2 \ot \ldots \ot 
h^{n-1}) \cr
&= & \ \sum \D^{n-1}  \wt S (h_{(1)}^1) \cdot h_{(2)}^1 \ot h^2 \ot 
\ldots \ot h^{n-1} \ot \sigma = \cr
&= & \ \sum S (h_{(1)(n)}^1)  h_{(2)}^1 \ot S (h_{(1)(n-1)}^1)  
h^2 \ot \ldots \cr
& &\ \ot S (h_{(1)(2)}^1)  h^{n-1} \ot \wt S (h_{(1)(1)}^1) \sigma \cr
&= & \ \sum \ve (h_{(n)}^1)  1 \ot S (h_{(n-1)}^1)  h^2 \ot \ldots 
\cr
& &\ \ot S (h_{(1)}^1)  h^{n-1} \ot \wt S (h_{(1)}^1) \sigma \cr
&= & \ 1 \ot S (h_{(n-1)}^1)  h^2 \ot \ldots \ot S (h_{(1)}^1)  
h^{n-1} \ot \wt S (h_{(1)}^1) \sigma \cr
&= & \ \delta_0  \tau_{n-1}  (h^1 \ot \ldots \ot h^{n-1})  , 
\end{eqnarray*}
and so forth.

Passing now to degeneracies,
\begin{eqnarray*}
&&\tau_n  \s_0 (h^1 \ot \ldots \ot h^{n+1}) = \ve (h^1)  \tau_n 
(h^2 \ot \ldots \ot h^{n+1}) = \cr
&&\ =\ \ve (h^1) \sum S (h_{(n)}^2)  h^3 \ot \ldots \ot S (h_{(2)}^2) 
 h^{n+1} \ot \wt S (h_{(1)}^2) \sigma  , 
\end{eqnarray*}
and on the other hand
\begin{eqnarray*}
& &\s_n  \tau_{n+1}^2 (h^1 \ot \ldots \ot h^{n+1}) = \cr
&= & \ \s_n \left( \sum S (h_{(n+1)}^2)  h^3 \ot \ldots \ot S 
(h_{(2)}^2) \sigma \ot \wt S (h_{(1)}^2)  \wt S^2 (h^1) \sigma \right) \cr
&= & \ \sum \ve (\wt S (h_{(1)}^2) \wt S^2 (h^1) \sigma )  S (h_{(n+1)}^2)  h^3 
\ot \ldots \ot S (h_{(2)}^2) \sigma \cr
&= & \ \ve (\sigma^{-1} \wt S^2 (h^1)\sigma) \sum \delta (h_{(1)}^2)  S (h_{(n+1)}^2)  h^3 \ot 
\ldots \ot S (h_{(2)}^2) \sigma \cr
&= & \ \ve (h^1)  S (h_{(n)}^2)  h^3 \ot \ldots \ot S (h_{(2)}^2) 
 h^{n+1} \ot \wt S (h_{(1)}^2) \sigma  . 
\end{eqnarray*}
In the next step
\begin{eqnarray*}
& &\ \tau_n  \s_1 (h^1 \ot \ldots \ot h^{n+1}) =  \ve (h^2)  
\tau_n (h^1 \ot h^3 \ot \ldots \ot h^{n+1}) \cr
&= & \ \ve (h^2) \cdot \D^{n-1}  \wt S (h^1) \cdot h^3 \ot \ldots 
\ot h^{n+1} \ot \sigma  , 
\end{eqnarray*}
while on the other hand
\begin{eqnarray*}
&& \ \s_0  \tau_{n+1} (h^1 \ot \ldots \ot h^{n+1}) = \cr
&& \ \sum \s_0 (S (h_{(n+1)}^1)  h^2 \ot \ldots \ot S (h_{(2)}^1)  
h^{n+1} \ot \wt S (h_{(1)}^1)) \sigma \cr
&= & \ \sum \ve (h^2) \cdot \ve (h_{(n+1)}^1) \cdot S (h_{(n)}^1)  h^3 
\ot \ldots \ot S (h_{(2)}^1)  h^{n+1} \ot \wt S (h_{(1)}^1) \sigma \cr
&= & \ \sum \ve (h^2) \cdot S (h_{(n-1)}^1)  h^3 \ot \ldots \ot S 
(h_{(2)}^1)  h^{n+1} \ot \wt S (h_{(1)}^1) \sigma  , 
\end{eqnarray*}
and similarly for $i = 2, \ldots n$. $\blacksquare$

The cohomology of the $(b,B)$-bicomplex corresponding to the cyclic 
module $\Hc_{(\delta,\sigma)}^{\natural}$ is, by definition, the 
\textit{cyclic cohomology}  $H  C_{(\delta,\sigma)}^* (\Hc)$ 
of $\Hc$ relative to the modular pair in involution ${(\delta,\sigma)}$.

\section{Examples}

The main point of the present paper is that thanks to $\sigma$ we remove 
the partial unimodularity condition of \cite{CM} 
on a Hopf algebra. 
As we shall see now, our general condition is fulfilled (modulo the passage to
a double cover) by the most popular Hopf algebras, 
including quantum groups and their duals.

\begin{proposition}
\noindent The following Hopf algebras are canonically endowed with 
modular pairs in involution: ribbon algebras, 
coribbon algebras and their tensor products, compact 
quantum groups in the sense of Woronowicz.
\end{proposition}

\noindent{\textit{Proof.}} If $\Hc$ is {\it quasitriangular}
with $R$-matrix $R$, then
$$
S^2 (h) = u \, h \, u^{-1} \, ,
$$
with
$$
u = \sum \, S ( R^{(2)}) \, R^{(1)} \ , \quad \varepsilon (u) = 1
$$
and
$$
\D \, u = (R_{21} \, R)^{-1} \, (u \otimes u) \, .
$$
\noindent By passage to a ``double cover'' \cite{RT}, i.e. an embedding in
$$
\Hc (\theta) = \Hc \, [\theta] / (\theta^2 - u \, S(u))
$$
one can assume that $u \, S(u) = S(u) \, u$, which is central, has a central 
square root $\theta$, such that
$$
\Delta (\theta) = (R_{21} \, R)^{-1} \, (\theta \otimes \theta) \ , \quad 
\varepsilon (\theta) = 1 \ , \quad S(\theta) = \theta \, .
$$
\noindent Taking
$$
\sigma = \theta^{-1} \, u \, ,
$$
one gets a canonical group-like element
$$
\Delta \, \sigma = \sigma \otimes \sigma \ , \quad \varepsilon (\sigma) = 1 \ , 
\quad S(\sigma) = \sigma^{-1} \, .
$$
\noindent  It is easy to check that the product of $\sigma^{-1}$ by the antipode
$$
 S' = \sigma^{-1} \cdot S \, ,
$$
satisfies the required condition 
$${S'}^2 = 1. $$
\noindent Indeed,
$$
\matrix{
{S'}^2 (h) &= &\sigma^{-1} \, S (\sigma^{-1} \, S (h)) = \sigma^{-1} \, S^2 (h) \, 
\sigma \hfill \cr
&= &\sigma^{-1} \, u \, h \, u^{-1} \, \sigma = \theta \, h \, \theta^{-1} = h 
\, . \hfill \cr
}
$$

\noindent This shows that $(\ve, \sigma)$ is a canonical modular pair in involution
for $\Hc$.
\smallskip

By dualizing the above definitions one obtains the notion of a
{\it coquasitriangular}, resp.  {\it coribbon} algebra. Among the most
prominent examples of coribbon algebras are the function
algebras of the classical quantum groups  $GL_q (N)$, $SL_q (N)$, $SO_q (N)$, 
$O_q (N)$ and $Sp_q (N)$.  
For a coribbon algebra $\Hc$, the analogue of
the above {\it ribbon group-like element} $\sigma$  is
the {\it ribbon character} $\delta \in \Hc^*$. The corresponding 
twisted antipode 
satisfies again the condition ${\widetilde S}^2 = 1$, so that $(\delta, 1)$
is a canonical modular pair in involution for $\Hc$.
\smallskip

Finally, for a compact quantum
group in the sense of \cite{W}, Theorem 5.6 of \cite{W} describing the
modular properties of the Haar measure shows that both the coordinate
algebra as well as its dual are provided with a canonical modular pair
in involution.
$\blacksquare$

\end{document}